\documentclass[12pt]{amsart}

\usepackage{latexsym}
\usepackage{amsthm,amsfonts,amssymb,amscd,mathrsfs}
\usepackage[leqno]{amsmath}
\usepackage{xspace}
\usepackage{cite}
\usepackage{fancyhdr}
\input{xypic}
\xyoption{all}

\newcommand{\numberseries}{\mdseries}   %Fontseries used for numbering theorem

\newlength{\thmtopspace}                %Space above theorem
\newlength{\thmbotspace}                %Space below theorem
\newlength{\thmheadspace}               %Space between theorem caption and text
\newlength{\thmindent}                  %For indenting

\newtheoremstyle{bfupright head,slanted body}
                {\thmtopspace}{\thmbotspace}
                {\slshape}{\thmindent}{\bfseries}{.}{\thmheadspace}
                {{\numberseries \thmnumber{(#2) }}\thmnote{#3}}

\newtheoremstyle{bfupright head,upright body}
                {\thmtopspace}{\thmbotspace}
                {\upshape}{\thmindent}{\bfseries}{.}{\thmheadspace}
                {{\numberseries \thmnumber{(#2) }}\thmnote{#3}}

\newtheoremstyle{bfit head,upright body}
                {\thmtopspace}{\thmbotspace}
                {\upshape}{\thmindent}{\upshape}{.}{\thmheadspace}
                {{\numberseries\thmnumber{(#2) }}
                {\bfseries\itshape\thmnote{\negthickspace#3}}}

\newtheoremstyle{it head,upright body}
                {\thmtopspace}{\thmbotspace}
                {\upshape}{\thmindent}{\upshape}{.}{\thmheadspace}
                {{\numberseries\thmnumber{(#2) }}
                {\itshape\thmnote{\negthickspace#3}}}

\newtheoremstyle{fixed bf head,slanted body}
                {\thmtopspace}{\thmbotspace}{\slshape}
                {\thmindent}{\bfseries}{.}{\thmheadspace}
                {{\numberseries \thmnumber{(#2) }}\thmname{#1}\thmnote{ (#3)}}

\newtheoremstyle{fixed bf head,upright body}
                {\thmtopspace}{\thmbotspace}{\upshape}
                {\thmindent}{\bfseries}{.}{\thmheadspace}
                {{\numberseries \thmnumber{(#2) }}\thmname{#1}\thmnote{ (#3)}}

\newtheoremstyle{indented paragraph}
                {\thmtopspace}{\thmbotspace}
                {\upshape}{\thmindent}{\upshape}{}{0pt}
                {\thmnote{#3 }}

\theoremstyle{bfupright head,slanted body}
\newtheorem{res}{}[section]             \newtheorem*{res*}{}

\theoremstyle{bfit head,upright body}
                 \newtheorem*{com*}{}

\theoremstyle{bfupright head,upright body}
\newtheorem{bfhpg}[res]{}               \newtheorem*{bfhpg*}{}

\theoremstyle{it head,upright body}
               \newtheorem*{ithpg*}{}

\theoremstyle{fixed bf head,slanted body}
\newtheorem{thm}[res]{Theorem}          \newtheorem*{thm*}{Theorem}
\newtheorem{prp}[res]{Proposition}      \newtheorem*{prp*}{Proposition}
\newtheorem{cor}[res]{Corollary}        \newtheorem*{cor*}{Corollary}
\newtheorem{lem}[res]{Lemma}            \newtheorem*{lem*}{Lemma}

\theoremstyle{fixed bf head,upright body}
\newtheorem{dfn}[res]{Definition}       \newtheorem*{dfn*}{Definition}
\newtheorem{obs}[res]{Observation}      \newtheorem*{obs*}{Observation}
\newtheorem{rmk}[res]{Remark}           \newtheorem*{rmk*}{Remark}
          \newtheorem*{exa*}{Example}
         \newtheorem*{exe*}{Exercise}
            \newtheorem{stp*}{Setup}
     \newtheorem*{dfns*}{Definitions}
    \newtheorem*{obss*}{Observations}
         \newtheorem*{rmks*}{Remarks}
        \newtheorem*{exas*}{Examples}

\theoremstyle{indented paragraph}

\newlength{\thmlistleft}        %leftmargin
\newlength{\thmlistright}       %rightmargin
\newlength{\thmlistpartopsep}   %partopsep
\newlength{\thmlisttopsep}      %topsep
\newlength{\thmlistparsep}      %parsep
\newlength{\thmlistitemsep}     %itemsep

\newcounter{eqc} 
                    {\end{list}}%

\newcounter{prt}
\newenvironment{prt}{\begin{list}{\upshape (\alph{prt})}%
                    {\usecounter{prt}%
                        \setlength{\leftmargin}{\thmlistleft}%
                        \setlength{\labelwidth}{\thmlistleft}%
                        \setlength{\rightmargin}{\thmlistright}%
                        \setlength{\partopsep}{\thmlistpartopsep}%
                        \setlength{\topsep}{\thmlisttopsep}%
                        \setlength{\parsep}{\thmlistparsep}%
                        \setlength{\itemsep}{\thmlistitemsep}}}%
                    {\end{list}}%

\newcounter{rqm}
\newenvironment{rqm}{\begin{list}{\upshape (\arabic{rqm})}%
                    {\usecounter{rqm}%
                        \setlength{\leftmargin}{\thmlistleft}%
                        \setlength{\labelwidth}{\thmlistleft}%
                        \setlength{\rightmargin}{\thmlistright}%
                        \setlength{\partopsep}{\thmlistpartopsep}%
                        \setlength{\topsep}{\thmlisttopsep}%
                        \setlength{\parsep}{\thmlistparsep}%
                        \setlength{\itemsep}{\thmlistitemsep}}}%
                    {\end{list}}%

                        {\end{list}}%

\newenvironment{itemlist}{\nopagebreak \begin{list}{$\bullet$}%
                       {\setlength{\leftmargin}{\thmlistleft}%
                        \setlength{\labelwidth}{\thmlistleft}%
                        \setlength{\rightmargin}{\thmlistright}%
                        \setlength{\partopsep}{\thmlistpartopsep}%
                        \setlength{\topsep}{\thmlisttopsep}%
                        \setlength{\parsep}{\thmlistparsep}%
                        \setlength{\itemsep}{\thmlistitemsep}}}%
                        {\end{list}}%

                        {\end{list}}%
\newlength{\myindent}
{\setlength{\myindent}{\parindent}\begin{list}{}%
                        {\setlength{\leftmargin}{#1}\setlength{\rightmargin}{#1}%
                        \setlength{\partopsep}{0pt}%
                        \setlength{\topsep}{\thmtopspace}%
                        \setlength{\parsep}{0pt}%
                        \setlength{\itemsep}{0pt}}
                        \item[]}
                        {\end{list}}%

\newenvironment{proof*}{\begin{proof}}{\renewcommand{\qed}{} \end{proof}}

\newcommand{\dispand}[1][and]{\hbox to \hsize{#1 \hfill} \nonumber \\}

\newlength{\seqsplit}
\title{Cotorsion pairs associated with Auslander categories}

\author{Edgar E. Enochs \ \ }

\address{\flushleft{Edgar} E. Enochs, Department of Mathematics, University of Kentucky, Lexington, KY, 40506-0027, USA} 

\email{enochs@ms.uky.edu} 

\urladdr{http://www.ms.uky.edu/\raisebox{-0.5ex}{\~{}}enochs/}

\author{\ \ Henrik Holm}

\address{\flushleft{Henrik} Holm, Department of Mathematical Sciences,
  University of Aarhus, Ny Munkegade, Building 530, DK--8000 Aarhus C,
  Denmark} 

\email{holm@imf.au.dk} 

\urladdr{http://home.imf.au.dk/holm/}

 

\newcommand{\A}{\mathcal{A}_C}
\newcommand{\B}{\mathcal{B}_C}

\newcommand{\Hom}{\operatorname{Hom}}
\newcommand{\Ext}{\operatorname{Ext}}
\newcommand{\Tor}{\operatorname{Tor}}
\newcommand{\Ker}{\operatorname{Ker}}
\renewcommand{\Im}{\operatorname{Im}}

\newcommand{\tensor}{\!\otimes\!}
\newcommand{\tensorR}{\!\otimes_R\!}

\newcommand{\Mod}{\mathsf{Mod}}
\newcommand{\Proj}{\mathcal{P}}
\newcommand{\Flat}{\mathcal{F}}
\newcommand{\Inj}{\mathcal{I}}

\begin{document}

\begin{abstract} 
  We prove that the Auslander class determined by a semidualizing
  module is the left half of a perfect cotorsion pair. We also prove
  that the Bass class determined by a semidualizing module is
  preenveloping.
\end{abstract}

\maketitle
\setcounter{section}{-1}

%-----------------------------------------------------------------------

\section{Introduction}

\noindent
The notion of semidualizing modules over commutative noetherian rings
goes back to Foxby \cite{HBFrelated} and Golod \cite{Golod}.
Christensen \cite{LWCsemi} extended this notion to semidualizing
complexes.

\medskip
\noindent
A semidualizing module or complex $C$ over a commutative
noe\-the\-rian ring gives rise to two full subcategories of the
derived category of the category of $R$--modules, namely the so-called
\textsl{Auslander class} and \textsl{Bass class} defined by
Avramov--Foxby \cite[(3.1)]{LLAHBF} and Christensen
\cite[def.\,(4.1)]{LWCsemi}.  Semidualizing complexes and their
Auslander/Bass classes have caught the attention of several authors,
see for example
\cite{LWCHH,EEEOMGJJZX,HBFrelated,LLAHBF,LWCsemi,HHPJ,EEESY,SSW,AFSSW,SISSW,LWCAFHH}.

\medskip
\noindent
Usually, one is interested in studying the \textsl{modules} in the
Auslander/Bass classes (by definition, the objects of these categories
are complexes), and in this paper we use $\A$ and $\B$ to denote the
ca\-te\-gories of all modules belonging to the Auslander class and
Bass class, respectively.

\medskip
\noindent
We mention that when $C$ itself is a (semidualizing) \textsl{module}
then one can describe $\A$ and $\B$ in terms of vanishing of certain
derived mo\-du\-le functors and invertibility of certain module
homomorphisms, see Avramov--Foxby \cite[prop.\,(3.4)]{LLAHBF} and
Christensen \cite[obs.\,(4.10)]{LWCsemi}.

\medskip
\noindent
In the case where $C$ is a \textsl{dualizing} module or complex, it is
possible to understand $\A$ and $\B$ in terms of the so-called
Gorenstein homological dimensions, see Enochs--Jenda--Xu
\cite{EEEOMGJJZX} and Christensen--Frankild--Holm \cite{LWCAFHH}.  A
similar description exists for other special semi\-dualizing complexes
$C$, see Holm--J{\o}rgensen \cite{HHPJ}.

\medskip
\noindent
In this paper we are concerned with what covering and enveloping
properties $\A$ and $\B$ possess. Our main results are Theorems
\eqref{thm:A} and \eqref{thm:B} which state, respectively, that:

\begin{res*}[Theorem A]
  Let $C$ be a semidualizing module over a commutative noetherian ring
  $R$. Then $(\A,(\A)^\perp)$ is a perfect cotorsion pair; in
  particular, $\A$ is covering. Furthermore, $\A$ is preenveloping.
\end{res*} 

\begin{res*}[Theorem B]
  Let $C$ be a semidualizing module over a commutative noetherian ring
  $R$. Then $\B$ is preenveloping.
\end{res*}

\noindent
As Corollary \eqref{cor:Gdim} we get:

\begin{res*} [Corollary C]
  Let $(R,\mathfrak{m},k)$ be a commutative, noetherian and local
  Co\-hen--Macaulay ring admitting a dualizing module. Then the
  following conclusions hold:
  \begin{prt}
  \item The class of all $R$--modules of finite Gorenstein
    projective/flat dimension is covering and preenveloping.
  \item The class of all $R$--modules of finite Gorenstein injective
    dimension is preenveloping.
  \end{prt}
\end{res*}

\section{Preliminaries} 

\noindent
In this section we briefly recall a number of definitions relevant to
this paper, namely the definitions of semidualizing modules, Auslander
categories, (pre)covers, (pre)envelopes, cotorsion pairs, and
Kaplansky classes. These notions will be used throughout the paper
without further mentioning.

\begin{bfhpg}[Setup]
  Throughout, $R$ is a fixed commutative noetherian ring with identity,
  and $C$ is a fixed semidualizing module for $R$, cf.~Definition
  \eqref{dfn:semi} below. We write $\Mod\,R$ for the category of
  $R$--modules.
\end{bfhpg}

\begin{rmk}
  Actually, we only need $R$ to be commutative and noetherian when we
  deal with semidualizing modules over $R$ and their Auslander/Bass
  classes. But in all of Section \ref{sec:exactsubcpx} for example,
  $R$ could be any ring.
\end{rmk}

\noindent
The next definition goes back to \cite{HBFrelated} (where the more
general \emph{PG--modules} are studied) and \cite{Golod}, but a more
recent reference is \cite[def.\,(2.1)]{LWCsemi}.

\begin{dfn} \label{dfn:semi}
  A \emph{semidualizing} module for $R$ is a finitely generated
  $R$--module $C$ such that:
  \begin{rqm}
  \item $\Ext^j_R(C,C)=0$ for all $j>0$, and
  \item The natural homothety morphism $\chi_C \colon R
    \longrightarrow \Hom_R(C,C)$ is an isomorphism.
  \end{rqm} 
\end{dfn} 

\noindent
The following lemma is straightforward.

\begin{lem} \label{lem:zero}
  Let $M$ be any $R$--module. If either $\Hom_R(C,M)=0$ or $C\tensorR
  M=0$, then $M=0$. \qed
\end{lem} 

\noindent
Next we recall the definitions, cf.~\cite[sec.\,1]{HBFrelated},
\cite[prop.\,(3.4)]{LLAHBF}, and \cite[def.\,(4.1) and
obs.\,(4.10)]{LWCsemi} of the ``module versions'' of the Auslander
categories with respect to the semidualizing module $C$.

\begin{dfn} \label{dfn:AB}
  The \emph{Auslander categories} \mbox{$\A=\A(R)$} and $\B=\B(R)$ are
  the full subcategories of $\Mod\,R$ whose objects are specified as
  follows:

  \medskip
  \noindent
  An $R$--module $M$ belongs to $\A$ if
  \begin{itemize}
  \item[{\small $(\mathrm{A}1)$}] $\Tor^R_j(C,M)=0$ for all $j>0$;
  \item[{\small $(\mathrm{A}2)$}] $\Ext^j_R(C,C \tensorR M)=0$ for all
    $j>0$;
  \item[{\small $(\mathrm{A}3)$}] $\eta_M \colon M \longrightarrow
    \Hom_R(C,C\tensorR M)$ is an isomorphism.
  \end{itemize}
  \smallskip
  \noindent  
  An $R$--module $N$ belongs to $\B$ if
  \begin{itemlist}
  \item[{\small $(\mathrm{B}1)$}] $\Ext^j_R(C,N)=0$ for all $j>0$;
  \item[{\small $(\mathrm{B}2)$}] $\Tor^R_j(C,\Hom_R(C,N))=0$ for all
    $j>0$;
  \item[{\small $(\mathrm{B}3)$}] $\varepsilon_N \colon C \tensorR
    \Hom_R(C,N) \longrightarrow N$ is an isomorphism.
  \end{itemlist}  
  \noindent 
  We often refer to $\A$ as the \emph{Auslander class}, and to $\B$ as
  the \emph{Bass class}.
\end{dfn}

\begin{rmk}
  In the notation of \cite{LLAHBF,LWCsemi}, the categories from
  Definition \eqref{dfn:AB} should have a subscript ``$0$''. However,
  since this paper only deals with modules (as opposed to complexes),
  and in order to keep notation as simple as possible, we have chosen
  to omit this ``$0$''.
\end{rmk}

\noindent
The next definition -- which is important for our main results
\eqref{thm:A}, \eqref{thm:B}, and \eqref{cor:Gdim} -- is taken
directly from Enochs--Jenda \cite[def.\,5.1.1 and 6.1.1]{EnochsRHA}.

\begin{dfn}
  Let $\mathcal{F}$ be a class of modules. An
  \emph{$\mathcal{F}$--precover} of a module $M$ is a homomorphism
  $\varphi \colon F \longrightarrow M$ with $F \in \mathcal{F}$ such
  that every homomorphism $\varphi' \colon F' \longrightarrow M$ with
  $F' \in \mathcal{F}$ factors as
  \begin{displaymath}
    \xymatrix{ {} & F' \ar[d]^-{\varphi'} \ar@{-->}[dl] \\
     F \ar[r]_-{\varphi} & M.
    }
  \end{displaymath}
  An $\mathcal{F}$--precover $\varphi \colon F \longrightarrow M$ is
  called an \emph{$\mathcal{F}$--cover} if every endomorphism $f
  \colon F \longrightarrow F$ with $\varphi f = \varphi$ is an
  automorphism. The class $\mathcal{F}$ is called \emph{(pre)covering}
  if every module has an $\mathcal{F}$--(pre)cover. There is a similar
  definition of \emph{$\mathcal{F}$--(pre)envelopes}.
\end{dfn}

\noindent
Our main Theorem \eqref{thm:A} uses the notion of a perfect cotorsion
pair. The definition, which is given below, is taken directly from
\cite[def.\,7.1.2]{EnochsRHA} and Enochs--L{\'o}pez-Ramos
\cite[def.~2.2]{EEEJALRkaplansky}.

\begin{dfn}
  A pair of module classes $(\mathcal{F},\mathcal{G})$ is a
  \emph{cotorsion pair} if $\mathcal{F}^\perp = \mathcal{G}$ and
  $\mathcal{F} = {}^\perp\mathcal{G}$, where
  \begin{align*}
    \mathcal{F}^\perp &= \{ N \in \mathsf{Mod}\,R \,|\, \Ext^1_R(F,N)=0
    \text{ for all } F \in \mathcal{F} \}, \text{ and} \\
    {}^\perp\mathcal{G} &= \{ M \in \mathsf{Mod}\,R \,|\, \Ext^1_R(M,G)=0
    \text{ for all } G \in \mathcal{G} \}.
  \end{align*}
  A cotorsion pair $(\mathcal{F},\mathcal{G})$ is called
  \emph{perfect} if $\mathcal{F}$ is covering and $\mathcal{G}$ is
  enveloping. 
\end{dfn}

\noindent
Finally, we need for Proposition \eqref{prp:Kaplansky} the notion of a
Kaplansky class. For the convenience of the reader we restate
\cite[def.~2.1]{EEEJALRkaplansky} below.

\begin{dfn} \label{dfn:Kaplansky} A class of modules $\mathcal{F}$ is
  called a \emph{Kaplansky class} if there exists a cardinal number
  $\lambda$ such that for every $x \in F \in \mathcal{F}$ there is a
  submodule $x \in F' \subseteq F$ with $|F'|\leqslant \lambda$ and
  $F', F/F' \in \mathcal{F}$.
\end{dfn}

\section{Exact subcomplexes of an exact complex} \label{sec:exactsubcpx}

\noindent
The main result of this section is Proposition
\eqref{prp:exactsubcomplexes}, which shows how to find desirable exact
subcomplexes of a given exact complex. This result is the cornerstone
in the proof of Proposition \eqref{prp:Kaplansky} in Section
\ref{sec:main}.

\begin{lem} \label{lem:cardinality} 
  Let $\lambda$ be any infinite cardinal number with $\lambda
  \geqslant |R|$, and let $M$ be any $R$--module. The following
  conclusions hold:
  \begin{prt}
  \item If $M$ is generated by a subset \mbox{$X\subseteq M$}
    satisfying $|X|\leqslant \lambda$ then also $|M|\leqslant
    \lambda$.
  \item If \mbox{$S\subseteq M$} is a submodule and with
    $|S|,|M/S|\leqslant \lambda$ then $|M|\leqslant \lambda$.
  \end{prt}
\end{lem}

\begin{proof}
  Part (a) is clear as there is an epimorphism \mbox{$R^{(X)}
    \longrightarrow M$}, and \mbox{$|R^{(X)}| \leqslant \lambda$}. For
  part (b) we pick a set of representatives \mbox{$X\subseteq M$} for
  the cosets of $S$ in $M$ (which has \mbox{$|X|=|M/S|\leqslant
    \lambda$}).  Then \mbox{$X \cup S \subseteq M$} generates $M$ and
  satisfies \mbox{$|X \cup S|\leqslant \lambda$}, so part (a) gives
  the desired conclusion.
\end{proof}

\begin{lem} \label{lem:functors}
  Let $\lambda$ be any infinite cardinal number. If $Q$ is a finitely
  generated module, and $P$ is a module with \mbox{$|P| \leqslant
    \lambda$} then the functors
  \begin{displaymath}
    \Hom_R(Q,-), P \tensorR- \colon \Mod\,R \longrightarrow \mathsf{Ab}
  \end{displaymath}
  have the property that the image of any module $A$ with
  \mbox{$|A|\leqslant \lambda$} has again car\-di\-na\-li\-ty
  $\leqslant \lambda$.
\end{lem}

\begin{proof}  
  Pick an integer $n>0$ and exact sequences
  \begin{displaymath}
    R^n \longrightarrow Q \longrightarrow 0 \quad \text{and} \quad 
    R^{(P)} \longrightarrow P \longrightarrow 0.
  \end{displaymath}
  Applying $\Hom_R(-,A)$ to the first of these sequences and
  $-\tensorR A$ to the second one, we get exactness of
  \begin{displaymath}
    0 \longrightarrow \Hom_R(Q,A) \longrightarrow \Hom_R(R^n,A) \cong
    A^n \ \text{ and}
  \end{displaymath}
  \begin{displaymath}
    A^{(P)} \cong R^{(P)}\tensorR A \longrightarrow P\tensorR A
    \longrightarrow 0.
  \end{displaymath}
  If \mbox{$|A|,|P| \leqslant \lambda$} then we have \mbox{$|A^n|,
    |A^{(P)}| \leqslant \lambda$}, and therefore also the desired
  conclusions, $|\Hom_R(Q,A)| \leqslant \lambda$ and $|P\tensorR
  A|\leqslant \lambda$.
\end{proof}

\begin{lem} \label{lem:shortcomplex}
  Let \mbox{$F \colon \Mod\,R \longrightarrow \mathsf{Ab}$} be an additive
  covariant functor which preserves direct limits, and assume that
  there exists an infinite cardinal number \mbox{$\lambda \geqslant
    |R|$} such that \mbox{$|FA| \leqslant \lambda$} for all modules
  $A$ with \mbox{$|A|\leqslant \lambda$}.  Let
  \begin{displaymath}
    \mathscr{E} = 
    E' \stackrel{d'}{\longrightarrow} E
    \stackrel{d}{\longrightarrow} E''
  \end{displaymath}
  be a complex of $R$--modules such that $F\mathscr{E}$ is exact.  Suppose
  \mbox{$S'\subseteq E'$}, \mbox{$S\subseteq E$}, \mbox{$S''\subseteq
    E''$} are submodules such that \mbox{$|S'|, |S|,|S''|\leqslant
    \lambda$}.  Then there exists a subcomplex,
  \begin{displaymath}
    \mathscr{T} =
    T' \stackrel{d'}{\longrightarrow} T
    \stackrel{d}{\longrightarrow} T''
  \end{displaymath}
  of $\mathscr{E}$ such that $F\mathscr{T}$ is exact, such that
  \mbox{$S'\subseteq T'$}, \mbox{$S\subseteq T$}, \mbox{$S''\subseteq
    T''$}, and such that $|T'|, |T|, |T''|\leqslant \lambda$.

  \medskip
  \noindent
  If, in addition, $\mathscr{E}$ is exact then we can choose
  $\mathscr{T}$ to be exact as well.
\end{lem}

\begin{proof}
  Replacing $S$ with \mbox{$S+d'(S')$} and $S''$ with
  \mbox{$S''+d(S)$} we see that we can assume that 
  \begin{displaymath}
    S' \stackrel{d'}{\longrightarrow} S
    \stackrel{d}{\longrightarrow} S''
  \end{displaymath}
  is a subcomplex of $\mathscr{E}$. Note that
  \begin{displaymath}
    \tag{\text{$\dagger$}}
    \mathscr{E} = 
    E' \stackrel{d'}{\longrightarrow} E
    \stackrel{d}{\longrightarrow} E'' = 
    \varinjlim \big( U' \stackrel{d'}{\longrightarrow} U
    \stackrel{d}{\longrightarrow} U'' \big)
  \end{displaymath}
  where the direct limit is taken over the family $\mathcal{U}$ of all
  subcomplexes \mbox{$U'\longrightarrow U\longrightarrow U''$} of
  $\mathscr{E}$ which contain \mbox{$S'\longrightarrow
    S\longrightarrow S''$} and satisfy that $U'/S'$, $U/S$ and
  $U''/S''$ are finitely generated. For each such subcomplex we also
  have $|U'|, |U|, |U''|\leqslant \lambda$ by Lemma
  \eqref{lem:cardinality}(b).

  \medskip
  \noindent
  Now, suppose that \mbox{$z\in \Ker(FS \longrightarrow FS'')
    \subseteq FS$}.  By $(\dagger)$ and the assumptions on the functor
  $F$ is follows that
  \begin{displaymath}
    \tag{\text{$\ddagger$}}
    FE' \longrightarrow FE \longrightarrow FE'' = 
    \varinjlim \big( FU' \longrightarrow FU \longrightarrow FU'' \big),
  \end{displaymath}
  which is exact by assumption. As the image of $z$ in $FE$
  belongs to
  \begin{displaymath}
    \Ker(FE \longrightarrow FE'') \,=\, \Im(FE'
    \longrightarrow FE), 
  \end{displaymath}
  the identity $(\ddagger)$ implies the existence of some
  \mbox{$U_z'\longrightarrow U_z \longrightarrow U_z''$} in
  $\mathcal{U}$ such that the image of $z$ in $FU_z$ belongs to
  \begin{displaymath} 
    \Im(FU_z' \longrightarrow FU_z).
  \end{displaymath}
  Then we define
  \begin{displaymath}
    T_0'\longrightarrow T_0\longrightarrow T_0'' = \sum_z
    \big(U_z' \longrightarrow U_z \longrightarrow U_z''\big) ,
  \end{displaymath}
  where the sum is taken over all \mbox{$z\in \Ker(FS \longrightarrow
    FS'')$}. By construction, there is an inclusion
  \begin{displaymath}
    \tag{\text{$\flat_0$}}
    \Im\! \big(\Ker(FS \longrightarrow FS'') \longrightarrow
    FT_0\big) \,\subseteq\, \Im(FT_0' \longrightarrow FT_0).
  \end{displaymath}
  We also note that the assumptions on $F$ imply that
  \begin{displaymath}
    \tag{\text{$\natural$}}
    |\Ker(FS \longrightarrow FS'')| \leqslant
    |FS| \leqslant \lambda,
  \end{displaymath}
  since $|S| \leqslant \lambda$. Consequently,
  \begin{displaymath}
    |T_0/S| = \Big|\sum_z U_z/S\Big| \leqslant \sum_z |U_z/S|
    \leqslant \sum_z \lambda \leqslant \lambda^2 = \lambda,
  \end{displaymath}
  where the penultimate inequality follows Lemma
  \eqref{lem:cardinality}(a), as $U_z/S$ is finitely generated, and
  the last inequality comes from $(\natural)$. Thus, we have
  \mbox{$|S|, |T_0/S| \leqslant \lambda$}, so Lemma
  \eqref{lem:cardinality}(b) implies that \mbox{$|T_0| \leqslant
    \lambda$}. Similarly,
  \begin{displaymath}
    |T_0'|, |T_0|, |T_0''| \leqslant \lambda.
  \end{displaymath}
  Now, going through the same procedure as above, but using the
  complex $T_0' \longrightarrow T_0\longrightarrow T_0''$ instead of
  $S'\longrightarrow S\longrightarrow S''$, we get a subcomplex
  \begin{displaymath}
    T_1'\longrightarrow T_1\longrightarrow T_1'' 
  \end{displaymath}
  of $\mathscr{E}$ containing $T_0' \longrightarrow T_0\longrightarrow
  T_0''$, and such that $|T_1'|, |T_1|, |T_1''| \leqslant \lambda$ and
  \begin{displaymath}
    \tag{\text{$\flat_1$}}
    \Im\! \big(\Ker(FT_0 \longrightarrow FT_0'') \longrightarrow
    FT_1\big) \,\subseteq\, \Im(FT_1' \longrightarrow FT_1).
  \end{displaymath}
  In this fashion we construct an increasing sequence
  \begin{displaymath}
    T_n'\longrightarrow T_n\longrightarrow T_n'' \quad , \quad n=0,1,2,\ldots
  \end{displaymath}
  of subcomplexes of $\mathscr{E}$ such that $|T_n'|, |T_n|, |T_n''|
  \leqslant \lambda$ and
  \begin{displaymath}
    \tag{\text{$\flat_n$}}
    \Im\! \big(\Ker(FT_{n-1} \longrightarrow FT_{n-1}'') \longrightarrow
    FT_n\big) \,\subseteq\, \Im(FT_n' \longrightarrow FT_n).
  \end{displaymath}
  Finally, we define a subcomplex of $\mathscr{E}$ by
  \begin{align*}
    T'\longrightarrow T\longrightarrow T'' &= \bigcup_{n=0}^{\infty}
    \big(T_n'\longrightarrow T_n \longrightarrow T_n''\big) \\ &=
    \varinjlim_{n\geqslant 0}
    \big(T_n'\longrightarrow T_n \longrightarrow T_n''\big).
  \end{align*}
  Note that \mbox{$|T| \leqslant \sum_{n \geqslant 0}|T_n| \leqslant
    \lambda + \lambda + \lambda+\cdots = \lambda$}, and similarly one
  gets \mbox{$|T'|, |T''| \leqslant \lambda$}. As $F$ commutes with
  direct limits we have
  \begin{displaymath}
    FT'\longrightarrow FT\longrightarrow FT'' =
    \varinjlim
    \big(FT_n'\longrightarrow FT_n \longrightarrow FT_n''\big).
  \end{displaymath}
  It is straightforward to verify that the conditions $(\flat_n)$
  ensure exactness of the complex above.
%  Assume that \mbox{$z \in \Ker(FT\longrightarrow FT'')$}. In
%  particular, \mbox{$z \in FT=\varinjlim FT_n$}, so there exist
%  \mbox{$m>0$} and \mbox{$z_m \in FT_m$} such that \mbox{$z_m
%    \longmapsto z$} under the map \mbox{$FT_m \longrightarrow FT$}.
%  Consider the image \mbox{$z_m'' \in FT_m''$} of $z_m$ under the map
%  \mbox{$FT_m \longrightarrow FT_m''$}. As \mbox{$z_m'' \longmapsto
%    0$} under the map \mbox{$FT_m'' \longrightarrow FT''=\varinjlim
%    FT_n''$} there is an \mbox{$\ell \geqslant m$} such that
%  \mbox{$z_m'' \longmapsto 0$} via \mbox{$FT_m'' \longrightarrow
%    FT_\ell''$}.  Then consider the image \mbox{$z_\ell \in FT_\ell$}
%  of $z_m$ under the map \mbox{$FT_m \longrightarrow FT_\ell$}, and
%  note that \mbox{$z_\ell \in \Ker(FT_\ell \longrightarrow
%    FT_\ell'')$}. Thus the image \mbox{$z_{\ell+1} \in FT_{\ell+1}$}
%  of $z_\ell$ under the map \mbox{$FT_\ell \longrightarrow
%    FT_{\ell+1}$} satisfies
%  \begin{displaymath}
%    z_{\ell+1} \in
%    \Im\! \big(\Ker(FT_{\ell} \longrightarrow FT_{\ell}'') \longrightarrow
%    FT_{\ell+1}\big),
%  \end{displaymath}
%  and thus by $(\flat_{\ell+1})$ there exists a \mbox{$z_{\ell+1}' \in
%    FT_{\ell+1}'$} such that \mbox{$z_{\ell+1}' \longmapsto
%    z_{\ell+1}$} via \mbox{$FT_{\ell+1}' \longrightarrow
%    FT_{\ell+1}$}. By construction, \mbox{$z_{\ell+1}' \longmapsto z$}
%  via \mbox{$FT_{\ell+1}' \longrightarrow FT$}, and hence $z$ belongs
%  to $\Im(FT'\longrightarrow FT)$.

   \medskip
   \noindent
   Concerning the last claim of the lemma we argue as follows: If, in
   addition, $\mathscr{E}$ is exact then we have exactness of
   $G\mathscr{E}$, where $G$ is the functor \mbox{$G=F \oplus
     \mathrm{id}$}, and $\mathrm{id} \colon \Mod\,R \longrightarrow
   \mathsf{Ab}$ is the forgetful functor. Hence, applying the part of
   the lemma which has already been established, but with $F$ replaced
   by $G$, we get that exactness of $G\mathscr{T}$.  Consequently,
   both $F\mathscr{T}$ and $\mathrm{id}\mathscr{T}=\mathscr{T}$ are
   exact.
\end{proof}

\begin{lem} \label{lem:pure}
  Let $\lambda$ be any infinite cardinal number with \mbox{$\lambda
    \geqslant |R|$}. If \mbox{$S \subseteq E$} is are $R$--modules
  such that \mbox{$|S|\leqslant \lambda$} then there is a pure
  submodule \mbox{$T\subseteq E$} of $E$ with \mbox{$S\subseteq T$}
  and \mbox{$|T|\leqslant \lambda$}.
\end{lem}

\begin{proof}
  A slight modification of the proof of \cite[Lemma 5.3.12]{EnochsRHA}
  gives this result.
\end{proof}

\begin{prp} \label{prp:exactsubcomplexes}
  Let \mbox{$F \colon \Mod\,R \longrightarrow \mathsf{Ab}$} be an additive
  covariant functor which preserves direct limits, and assume that
  there exists an infinite cardinal number \mbox{$\lambda \geqslant
    |R|$} such that \mbox{$|FA| \leqslant \lambda$} for all modules
  $A$ with \mbox{$|A|\leqslant \lambda$}.  Let
  \begin{displaymath}
    E = \cdots \longrightarrow E_{n+1} \longrightarrow E_{n}
    \longrightarrow E_{n-1} \longrightarrow \cdots 
  \end{displaymath}
  be an exact complex of $R$--modules such that also $FE$ is exact.
  If \mbox{$S_n\subseteq E_n$} is a submodule such that
  \mbox{$|S_n|\leqslant \lambda$} for each $n \in \mathbb{Z}$, then
  there is an exact subcomplex
  \begin{displaymath}
    D = \cdots \longrightarrow D_{n+1} \longrightarrow D_{n}
    \longrightarrow D_{n-1} \longrightarrow \cdots 
  \end{displaymath}
  of $E$ such that also $FD$ is exact. Furthermore, $D_n\subseteq E_n$
  is a pure submodule, $S_n\subseteq D_n$, and \mbox{$|D_n|\leqslant
    \lambda$} for each $n \in \mathbb{Z}$.
\end{prp}

\begin{proof}
  For each $n$ we construct a chain \mbox{$S_n=S_n^0\subseteq
    S_n^1\subseteq S_n^2\subseteq \cdots$} of submodules of $E_n$ with
  $|S_n^i|\leqslant \lambda$ for all $n$ and $i$ as follows:

  \medskip
  \noindent
  First pick a function \mbox{$f \colon \mathbb{N}_0 \longrightarrow
    \mathbb{Z}$} with the property that for each \mbox{$n \in
    \mathbb{Z}$} the set \mbox{$f^{-1}(\{n\})=\{k \in \mathbb{N}_0
    \,|\, f(k)=n\}$} is infinite.  Of course, we are given that
  $S_n^0=S_n$ for all $n$.
  
  \medskip
  \noindent
  Assume that we for some $k\geqslant 0$ have constructed $S_n^{2k}$
  for all $n$ with the property that \mbox{$|S_n^{2k}|\leqslant
    \lambda$}.  By Lemma \eqref{lem:pure} we get a pure submodule
  \mbox{$S_n^{2k+1} \subseteq E_n$} which contains $S_n^{2k}$ and
  which satisfies \mbox{$|S_n^{2k+1}|\leqslant \lambda$}.

  \medskip
  \noindent
  Now assume that we for some $k\geqslant 0$ have constructed
  $S_n^{2k+1}$ for all $n$ with the property that
  \mbox{$|S_n^{2k+1}|\leqslant \lambda$}. Then we construct
  $S_n^{2k+2}$ by the following procedure: Consider the complex
  \begin{displaymath}
    \mathscr{E}_k = 
    E_{f(k)+1} \longrightarrow E_{f(k)} \longrightarrow E_{f(k)-1}.
  \end{displaymath}
  By hypothesis both $\mathscr{E}_k$ and $F\mathscr{E}_k$ are exact.
  Applying Lemma \eqref{lem:shortcomplex} to this complex
  $\mathscr{E}_k$ and to the submodules,
  \begin{displaymath}
    S_{f(k)+1}^{2k+1} \subseteq E_{f(k)+1} \ , \ \ 
    S_{f(k)}^{2k+1} \subseteq E_{f(k)} \ \text{ and } \ \ 
    S_{f(k)-1}^{2k+1} \subseteq E_{f(k)-1}
  \end{displaymath}
  we get an exact subcomplex
  \begin{displaymath}
    \mathscr{T}_k =
    T' \longrightarrow T \longrightarrow T''
  \end{displaymath}
  of $\mathscr{E}_k$ where also $F\mathscr{T}_k$ is exact, and furthermore,
  \begin{displaymath}
    S_{f(k)+1}^{2k+1} \subseteq T' \ , \ \ 
    S_{f(k)}^{2k+1} \subseteq T \ \text{ and } \ \ 
    S_{f(k)-1}^{2k+1} \subseteq T'',
  \end{displaymath}
  and $|T'|, |T|, |T''| \leqslant \lambda$. We then define
  \begin{displaymath}
    S_{f(k)+1}^{2k+2}=T' \ , \ \ 
    S_{f(k)}^{2k+2}=T \ \text{ and } \ \ 
    S_{f(k)-1}^{2k+2}=T'',
  \end{displaymath}
  and let $S_n^{2k+2}=S_n^{2k+1}$ for all \mbox{$n \notin \{f(k)+1,
    f(k), f(k)-1\}$}.
  
  \medskip 
  \noindent
  Having constructed the sequences \mbox{$S_n=S_n^0\subseteq
    S_n^1\subseteq S_n^2\subseteq \cdots \subseteq E_n$} for $n \in
  \mathbb{Z}$ as above, we set \mbox{$D_n =\bigcup_{i=0}^{\infty}
    S_n^i$}. Clearly, \mbox{$S_n\subseteq D_n$} and
  \mbox{$|D_n|\leqslant \lambda$}.  Since
  \mbox{$D_n=\bigcup_{k=0}^{\infty} S_n^{2k+1}$} and since
  \mbox{$S_n^{2k+1}\subseteq E_n$} is a pure submodule for each $k$,
  it follows that \mbox{$D_n \subseteq E_n$} is a pure submodule.

  \medskip
  \noindent
  For each $n$, the differential \mbox{$E_n \longrightarrow E_{n-1}$}
  restricts to a homomorphism \mbox{$D_n \longrightarrow D_{n-1}$}: If
  $x \in D_n$ then there exists an $i_0$ such that $x \in S_n^i$ for
  all $i \geqslant i_0$. Since $f^{-1}(\{n\})$ is infinite, there
  exists $k \geqslant 0$ satisfying both $f(k)=n$ and $2k+1 \geqslant
  i_0$.  By our construction, 
  \begin{displaymath}
    S_{n+1}^{2k+2} \longrightarrow S_n^{2k+2} \longrightarrow
    S_{n-1}^{2k+2} 
  \end{displaymath}
  is a subcomplex of \mbox{$E_{n+1} \longrightarrow E_n
    \longrightarrow E_{n-1}$}, so in particular the diffe\-rential
  \mbox{$E_n \longrightarrow E_{n-1}$} maps \mbox{$x \in S_n^{2k+1}
    \subseteq S_n^{2k+2}$} into \mbox{$S_{n-1}^{2k+2} \subseteq
    D_{n-1}$}. Hence
  \begin{displaymath}
    D = \cdots \longrightarrow D_{n+1} \longrightarrow D_{n}
    \longrightarrow D_{n-1} \longrightarrow \cdots 
  \end{displaymath}
  is a subcomplex of $E$. In fact, for the $n$'th segment of $D$ we
  have the expression
  \begin{displaymath}
    D_{n+1} \longrightarrow D_{n} \longrightarrow D_{n-1} =
    \varinjlim_{k \in f^{-1}(\{n\})} \big( S_{n+1}^{2k+2}
    \longrightarrow S_n^{2k+2} \longrightarrow S_{n-1}^{2k+2} \big),
  \end{displaymath}
  and since each of the complexes 
  \begin{displaymath}
     S_{n+1}^{2k+2} \longrightarrow S_n^{2k+2} \longrightarrow
     S_{n-1}^{2k+2} \quad , \quad k \in f^{-1}(\{n\})
  \end{displaymath}
  are exact and stay exact when we apply $F$ to them, the same is true
  for $D_{n+1} \longrightarrow D_{n} \longrightarrow D_{n-1}$, as $F$
  commutes with direct limits.
\end{proof}

\section{Covers and envelopes by Auslander categories} \label{sec:main}

\noindent
This last section is concerned with covering and enveloping properties
of the Auslander categories. Our main results are Theorems
\eqref{thm:A} and \eqref{thm:B}.

\medskip
\noindent
To prove our main theorems, we need the alternative descriptions of
the modules in the Auslander categories given in Propositions
\eqref{prp:characterizingA} and \eqref{prp:characterizingB}. To this
end, we introduce two new classes of modules:

\begin{dfn}
  The classes of \emph{$C$--injective} and \emph{$C$--flat} modules
  are defined as
  \begin{align*}
    &\Inj_C = \Inj_C(R) = \{\Hom_R(C,I) \,|\, \text{$I$ injective
      $R$--module}\}, \\ 
    &\Flat_C = \Flat_C(R) = \{C\tensorR F \,|\, \text{$F$ flat
      $R$--module}\}. 
  \end{align*}
\end{dfn} 

\begin{obs} 
  $R$ is a semidualizing module for itself, and by setting $C=R$ in
  the de\-fi\-nition above we see that \mbox{$\Inj_R$} and $\Flat_R$
  are the classes of (ordinary) injective and flat $R$--modules,
  respectively.
\end{obs}

\noindent
The proof of the next lemma is straightforward.

\begin{lem} \label{lem:description}
  For $R$--modules $U$ and $V$ one has the implications:
  \begin{prt}
  \item $U \in \Inj_C \iff U \in \A$ and $C\tensorR U \in \Inj_R$.
  \item $V \in \Flat_C \iff V \in \B$ and $\Hom_R(C,V) \in
    \Flat_R$. \qed
  \end{prt}
\end{lem}

\begin{rmk} \label{rmk:Proj}
  The classes $\Inj_C$, $\Flat_C$ and also
  \begin{align*}
    \Proj_C = \Proj_C(R) = \{C\tensorR P \,|\, P \text{ projective
      module}\} 
  \end{align*}
  were used in \cite{HHPJ}, and it was proved in
  \cite[lem.\,2.14]{HHPJ} (compared with
  \cite[thm.\,2.5]{EEEJALRkaplansky}) that $\Flat_C$ is preenveloping.
  When $R$ is a Cohen--Macaulay local ring and $C$ is a dualizing
  module for $R$ it was proved in \cite[prop.\,1.5]{EEEOMGJJZX} that
  $\Inj_C$ is preenveloping and $\Proj_C$ is precovering.
\end{rmk}

\noindent
We will make use of the following:

\begin{prp} \label{prp:envcov} 
  $\Inj_C$ is enveloping and $\Flat_C$ is covering. In particular, for
  any $R$--module $M$ there exist complexes
  \begin{displaymath}
    \mathscr{U} = 
    0 \longrightarrow M \longrightarrow U^0 \longrightarrow U^1
    \longrightarrow U^2 \longrightarrow \cdots 
  \end{displaymath}
  with $U^0,U^1,U^2,\ldots \in \Inj_C$, and 
  \begin{displaymath}
    \mathscr{V} = 
    \cdots \longrightarrow V_2 \longrightarrow V_1 \longrightarrow V_0
    \longrightarrow M \longrightarrow 0 
  \end{displaymath}
  with \mbox{$V_0,V_1,V_2 \in \Flat_C$} such that
  \mbox{$C\tensorR \mathscr{U}$} and \mbox{$\Hom_R(C,\mathscr{V})$} are exact.
\end{prp} 

\begin{proof}
  Even if $\Inj_C$ is just preenveloping and $\Flat_C$ is just
  precovering there will exist complexes $\mathscr{U}$ and $\mathscr{V}$ with
  \mbox{$U^i \in \Inj_C$} and \mbox{$V_j \in \Flat_C$} such that
  \begin{eqnarray*}
    &\Hom(\mathscr{U},\Hom(C,I)) \cong \Hom(C\tensor \mathscr{U},I),
    \text{ and} 
    \\
    &\Hom(C\tensor F,\mathscr{V})
  \end{eqnarray*}
  are exact for all injective modules $I$ and all flat modules $F$.
  Taking $I$ to be faithfully injective and \mbox{$F=R$}, we see that
  \mbox{$C\tensor \mathscr{U}$} and \mbox{$\Hom(C,\mathscr{V})$} are exact.

  \medskip
  \noindent
  Thus, the proposition is proved when we have argued that $\Inj_C$
  is enveloping and $\Flat_C$ is covering. The class $\Inj_R$ is
  known to be enveloping by Xu \cite[thm.\,1.2.11]{Xu}, since injective
  hulls in the sense of Eckmann--Schopf \cite{EckmannSchopf} always
  exists. The class $\Flat_R$ is known to be covering by
  Bican--Bashir--Enochs \cite{BicanBashirEnochs}. Now, it is easy to
  see that for any module $M$, the composition
  \begin{displaymath}
    M \longrightarrow \Hom(C,C\tensor M) \longrightarrow
    \Hom(C,E(C\tensor M)) 
  \end{displaymath}
  is an $\Inj_C$--envelope of $M$, where $E(-)$ denotes the injective
  envelope. Likewise, the composition
  \begin{displaymath}
    C\tensor F(\Hom(C,M)) \longrightarrow C\tensor \Hom(C,M)
    \longrightarrow M
  \end{displaymath}
  is a $\Flat_C$--cover of $M$, where $F(-)$ denotes the flat cover.
\end{proof} 

\noindent
Having established Lemma \eqref{lem:description} and Proposition
\eqref{prp:envcov}, the proof of the next result is similar to that of
\cite[prop.\,5.5.4]{Xu}.

\begin{prp} \label{prp:characterizingA}
  A module $M$ belongs to $\A$ if and only if there exists an exact
  sequence
  \begin{displaymath}
    \tag{\text{$\dagger$}}
    \cdots \longrightarrow F_2 \longrightarrow F_1 \longrightarrow F_0
    \longrightarrow U^0 \longrightarrow U^1 \longrightarrow U^2
    \longrightarrow \cdots  
  \end{displaymath}
  satisfying the following conditions:
  \begin{rqm}
  \item \mbox{$F_0, F_1, F_2,\ldots \in \Flat_R$} and \mbox{$U^0,
      U^1, U^2, \ldots \in \Inj_C$};
  \item $M = \Ker(U^0 \longrightarrow U^1)$;
  \item \mbox{$C\tensorR (\dagger)$} is exact. \qed
  \end{rqm}
\end{prp}

\noindent
Dually one proves the next result which is similar to
\cite[prop.\,5.5.5]{Xu}.

\begin{prp} \label{prp:characterizingB}
  A module $M$ belongs to $\B$ if and only if there exists an exact
  sequence
  \begin{displaymath}
    \tag{\text{$\ddagger$}}
    \cdots \longrightarrow V_2 \longrightarrow V_1 \longrightarrow V_0
    \longrightarrow I^0 \longrightarrow I^1 \longrightarrow I^2
    \longrightarrow \cdots  
  \end{displaymath}
  satisfying the following conditions:
  \begin{rqm}
  \item \mbox{$V_0, V_1, V_2,\ldots \in \Flat_C$} and \mbox{$I^0,
      I^1, I^2, \ldots \in \Inj_R$};
  \item $M = \Ker(I^0 \longrightarrow I^1)$;
  \item \mbox{$\Hom_R(C,(\ddagger))$} is exact. \qed
  \end{rqm}
\end{prp}

\begin{rmk}
  Taking the necessary precautions, one can study Auslander categories
  over non-noetherian rings. In this generality one can also prove
  versions of for example Propositions \eqref{prp:characterizingA} and
  \eqref{prp:characterizingB}, see \cite{HHDW}.
\end{rmk}

\medskip
\noindent
In addition to the fact that \cite[chap.\,5.5]{Xu} assumes $C$ to be
dualizing (and not just semidualizing), there is another important
difference between \cite[prop.\,5.5.5]{Xu} and Proposition
\eqref{prp:characterizingB}: Namely, we work with $\Flat_C$ whereas Xu
works with $\Proj_C$ (which he denotes $\mathcal{W}$); see Remark
\eqref{rmk:Proj}.

\medskip
\noindent
From our point of view $\Flat_C$ is more flexible than $\Proj_C$.  For
example, $\Flat_C$ is closed under pure submodules and pure quotients;
see Proposition \eqref{prp:pure} below, whereas $\Proj_C$ in general
does not have these properties.

%\begin{proof}
%  Similar to the proof of Proposition \eqref{prp:characterizingA}.
%\end{proof}

\begin{prp} \label{prp:pure}
  The classes $\Inj_C$ and $\Flat_C$ are closed under pure
  submodules and pure quotients.
\end{prp}

\begin{proof}
  First we prove that $\Inj_C$ is closed under pure submodules and
  pure quotients. To this end let
  \begin{displaymath}
    0 \longrightarrow M \longrightarrow \Hom(C,I) \longrightarrow N
    \longrightarrow 0 
  \end{displaymath}
  be a pure exact sequence with $I$ injective. Applying
  \mbox{$C\tensor -$} to this sequence we get another pure exact
  sequence,
  \begin{displaymath}
    \tag{\text{$*$}}
    0 \longrightarrow C\tensor M \longrightarrow I \longrightarrow
    C\tensor N \longrightarrow 0.
  \end{displaymath}
  As $I$ is injective, and as $\Inj_R$ is closed under pure submodules
  and pure quotients, we conclude that the modules
  \mbox{$C\tensor M$} and \mbox{$C\tensor N$} are
  injective. Applying $\Hom(C,-)$ to the pure exact sequence $(*)$
  we get exactness of the lower row in the following commutative
  diagram:
  {\small
  \begin{displaymath}
    \xymatrix{0 \ar[r] & M \ar[d]^-{\eta_M} \ar[r] & \Hom(C,I)
      \ar@{=}[d] \ar[r] & 
      N \ar[d]^-{\eta_N} \ar[r] & 0 \\
    0 \ar[r] & \Hom(C,C\tensor  M) \ar[r] & \Hom(C,I) \ar[r] &
      \Hom(C,C\tensor  N) \ar[r] & 0 }
  \end{displaymath}
}Since \mbox{$C\tensor M$} and \mbox{$C\tensor N$} are injective by
Lemma \eqref{lem:description}(a), we are done if we can prove that
$\eta_M$ and $\eta_N$ are isomorphisms.  By the snake-lemma, $\eta_M$
is injective, $\eta_N$ is surjective, and \mbox{$\operatorname{Ker}
  \eta_N \cong \operatorname{Coker} \eta_M$}.  Hence, it suffices to
argue that $\operatorname{Coker} \eta_M = 0$, and by Lemma
\eqref{lem:zero} it is enough to show that 
  \begin{displaymath}
    C\otimes \operatorname{Coker} \eta_M =0.
  \end{displaymath}
  Right exactness of \mbox{$C\tensor -$} gives exactness of
  \begin{displaymath}
    \xymatrix{ C\tensor M \ar[r]^-{C\otimes  \eta_M} &
      C\tensor  \Hom(C,C\tensor  M) \ar[r] & C\tensor  \operatorname{Coker}
      \eta_M \ar[r] & 0, 
    }
  \end{displaymath}
  Since \mbox{$C\tensor M$} is injective, $C\tensor \eta_M$ is an
  isomorphism with inverse $\varepsilon_{(C\tensor M)}$. In
  particular, $C\tensor \eta_M$ is surjective, and hence $C\tensor
  \operatorname{Coker} \eta_M =0$, as desired.

  \medskip
  \noindent
  Similarly, as the class of flat modules is closed under pure
  submodules and pure quotients, one proves that $\Flat_C$ also has
  these properties.
\end{proof}

\begin{prp} \label{prp:Kaplansky}
  $\A$ and $\B$ are Kaplansky classes.
\end{prp}

\begin{proof}
  We only prove that $\A$ is Kaplansky, as the proof for $\B$ is
  similar. We claim that any infinite cardinal number $\lambda
  \geqslant |R|$ implements the Kaplansky property for $\A$,
  cf.~Definition \eqref{dfn:Kaplansky}:

  \medskip
  \noindent
  Let \mbox{$x \in M \in \A$}. By Proposition
  \eqref{prp:characterizingA} there is an exact sequence
  \begin{displaymath}
    E =
    \cdots \longrightarrow F_2 \longrightarrow F_1 \longrightarrow F_0
    \longrightarrow U^0 \longrightarrow U^1 \longrightarrow U^2
    \longrightarrow \cdots  
  \end{displaymath}
  satisfying the conditions (1), (2), and (3) of that result. Now
  consider the submodules:
  \begin{itemlist}
  \item $Rx\subseteq U^0$ (which has $|Rx| \leqslant \lambda$);
  \item $0\subseteq U^n$ for $n\geqslant 1$;
  \item $0\subseteq F_n$ for $n\geqslant 0$.
  \end{itemlist}
  Applying Proposition \eqref{prp:exactsubcomplexes} with
  \mbox{$F=C\tensor-$} (cf.~Lemma \eqref{lem:functors}) to this
  situation we get an exact subcomplex
  \begin{displaymath}
    D =
    \cdots \longrightarrow G_2 \longrightarrow
    G_1 \longrightarrow G_0 
    \longrightarrow W^0 \longrightarrow W^1
    \longrightarrow W^2 
    \longrightarrow \cdots  
  \end{displaymath}
  of $E$ with \mbox{$Rx \subseteq W^0$} and such that \mbox{$C\tensor
    D$} is exact, and furthermore, \mbox{$G_n \subseteq F_n$} and
  \mbox{$W^n \subseteq U^n$} are pure submodules, and \mbox{$|G_n|,
    |W^n| \leqslant \lambda$}.

  \medskip
  \noindent
  By the condition \eqref{prp:characterizingA}(1), $F_n \in \Flat_R$
  and $U^n \in \Inj_C$, and thus Proposition \eqref{prp:pure}
  implies that 
  \begin{displaymath}
    \tag{\text{$*$}}
    G_n, F_n/G_n \in \Flat_R
    \quad \text{and} \quad
    W^n, U^n/W^n \in \Inj_C.
  \end{displaymath}
  Hence, Proposition \eqref{prp:characterizingA} implies that
  \mbox{$M'=\Ker(W^0 \longrightarrow W^1)$} belongs to $\A$.  As
  \mbox{$x \in M=\Ker(U^0 \longrightarrow U^1)$} we also have \mbox{$x
    \in M'$}, and of course \mbox{$|M'| \leqslant |W^0| \leqslant
    \lambda$}. Thus, it remains to argue that $M/M' \in \A$. By
  construction we have a pure exact sequence of complexes
  \begin{displaymath}
    \tag{\text{$\sharp$}}
    0 \longrightarrow D \longrightarrow E
    \longrightarrow E/D \longrightarrow 0
  \end{displaymath}
  As $D$ and $E$ are exact then so is $E/D$, in particular it follows
  that \mbox{$M/M' \cong \Ker(U^0/W^0 \longrightarrow U^1/W^1)$}.
  Purity of $(\sharp)$ gives exactness of
  \begin{displaymath}
    0 \longrightarrow C\tensor D \longrightarrow C\tensor E
    \longrightarrow C\tensor(E/D) \longrightarrow 0,
  \end{displaymath}
  and since $C \tensor D$ and $C\tensor E$ are exact then so is
  $C\tensor(E/D)$. By $(*)$ the complex $E/D$ consists of modules of
  the form \eqref{prp:characterizingA}(1), and these arguments show
  that $M/M' \in \A$, as desired.
\end{proof}

%\begin{proof}[Sketch of a possible proof?]
%  \textsl{The idea which hopefully works is as follows:} By
%  Propositions \eqref{prp:characterizingA} and
%  \eqref{prp:characterizingB} and by Proposition \eqref{prp:pure}, the
%  modules in $\A$ and $\B$ can be described in terms of certain exact
%  sequences consisting of modules from certain classes which are
%  closed under pure submodules and pure quotients. Now one can apply
%  the ``zig-zag'' procedure from the proof of Enochs--L{\'o}pez-Ramos
%  \cite[prop.\,2.6]{EEEJALRkaplansky}
%\end{proof}

\begin{thm} \label{thm:A} 
  $(\A,(\A)^\perp)$ is a perfect cotorsion pair; in particular, $\A$
  is covering. Furthermore, $\A$ is preenveloping.
\end{thm} 

\begin{proof}
  By Proposition \eqref{prp:Kaplansky}, the class $\A$ is Kaplansky.
  Clearly, $\A$ contains the projective modules, and is closed under
  extensions and direct limits. Hence
  \cite[thm.\,2.9]{EEEJALRkaplansky} implies that $(\A,(\A)^\perp)$ is
  a perfect cotorsion pair. As $\A$ is also closed under products,
  \cite[thm.\,2.5]{EEEJALRkaplansky} gives that $\A$ is preenveloping.
\end{proof}

\begin{thm} \label{thm:B}
  $\B$ is preenveloping.
\end{thm}

\begin{proof}
  By Proposition \eqref{prp:Kaplansky}, the class $\B$ is Kaplansky.
  Since $\B$ is closed under direct limits and products,
  \cite[thm.\,2.5]{EEEJALRkaplansky} implies that $\B$ is
  preenveloping.
\end{proof}

\begin{cor} \label{cor:Gdim}
  Let $(R,\mathfrak{m},k)$ be a commutative, noetherian and local
  Cohen--Macaulay ring admitting a dualizing module. Then the
  following conclusions hold:
  \begin{prt}
  \item The class of $R$--modules of finite Gorenstein projective/flat
    dimension is covering and preenveloping.
  \item The class of $R$--modules of finite Gorenstein injective
    dimension is preenveloping.
  \end{prt}
\end{cor}

\begin{proof}
  Taking $C$ to be the dualizing module for $R$, the assertions follow
  immediately from comparing Theorems \eqref{thm:A} and \eqref{thm:B}
  with \cite[cor.\,2.4 and 2.6]{EEEOMGJJZX}.
\end{proof}

\section*{Acknowledgements}

\noindent
We thank Diana White for useful comments, and for correcting many
typos in an ealier version of this manuscript.

\providecommand{\bysame}{\leavevmode\hbox to3em{\hrulefill}\thinspace}
\providecommand{\MR}{\relax\ifhmode\unskip\space\fi MR }
\providecommand{\MRhref}[2]{%
  \href{http://www.ams.org/mathscinet-getitem?mr=#1}{#2}
}
\providecommand{\href}[2]{#2}

\bigskip
\bigskip

\end{document}